\documentclass[12pt,draft]{article}

\usepackage{cite}
\usepackage[cp1251]{inputenc}
\usepackage[T2A]{fontenc}
\usepackage[english,russian,ukrainian]{babel}
\usepackage{amsmath,amsfonts,amssymb}
\usepackage{geometry}
\begin{document}

{\large\textbf{Умови регулярності розв'язків деяких параболічних систем}}

\medskip

{\large\textbf{Regularity conditions for solutions to some parabolic systems}}

\medskip

\begin{flushleft}
\small{\textbf{О. В. Дяченко (O. V. Diachenko)}}

\small{(Національний технічний університет України "'Київський політехнічний інститут імені Ігоря Сікорського"')}
\small{ol$\underline{\phantom{k}}$v$\underline{\phantom{k}}$dyachenko@ukr.net}

\medskip

\small{\textbf{В. М. Лось (V. M. Los)}} 

\small{(Національний технічний університет України "'Київський політехнічний інститут імені Ігоря Сікорського"')}
\small{v$\underline{\phantom{k}}$\,los@yahoo.com}

\end{flushleft}

\normalsize

Досліджено глобальну і локальну регулярність узагальнених розв'язків
початково-крайової задачі для параболічної за Петровським системи диференціальних рівнянь другого порядку.
Результати сформульовано в термінах приналежності правих частин задачі деяким узагальненим просторам Соболєва.
Отримано нові достатні умови класичності узагальненого розв'язку.

\smallskip

We investigate global and local regularity of generalized solutions to
parabolic initial-boundary value problem for Petrovskii system of second order differential equations.
Results are formulated in terms of the belonging of right-hand sides of the problem to some generalized Sobolev spaces.
We also obtain new sufficient conditions under which the generalized solution should be classical.

\bigskip

\textbf{1. Вступ.}
Мета цієї праці - доповнити результати нашої статті \cite{DiachenkoLos22JEPE} про коректну розв'язність деяких параболічних початково-крайових задач
в узагальнених функціональних просторах Соболєва теоремами про достатні умови регулярності розв'язків задач. Ми досліджуємо параболічні за Петровським системи
диференціальних рівнянь другого порядку, задані у багатовимірному скінченному циліндрі $\Omega$ з гладкою бічною поверхнею.
Регулярність розв'язків цих систем характеризуємо у термінах гільбертових анізотропних функціональних просторів $H^{s,s/2;\varphi}(\Omega)$, де дійсне число
$s\geq2$, а функціональний параметр $\varphi:[1,\infty)\rightarrow(0,\infty)$ повільно змінюється на нескінченності у сенсі Карамати.
Параметр $\varphi$ уточнює основну степеневу регулярність розв'язків, яка задається числовими параметрами $s$ і $s/2$ за просторовими і часовою змінними відповідно.

Вказані простори є окремим випадком просторів, уведених Л. Хермандером \cite{Hermander63,Hermander83}. Якщо $\varphi(\cdot)\equiv1$, то $H^{s,s/2;\varphi}(\Omega)$ --
анізотропні простори Соболєва $H^{s,s/2}(\Omega)$, які разом з просторами Гельдера широко застосовуються у теорії параболічних задач
\cite{Solonnikov65,LionsMagenes72ii, Ivasyshen90, ZhitarashuEidelman98,Ilin60, IlinKalashOleinik62}.
Її центральний результат стверджує, що ці задачі коректно розв'язні (за Адамаром) на парах відповідних анізотропних просторів Соболєва і Гельдера,
тобто породжують ізоморфізми на цих парах.

У теперішній час параболічні задачі активно досліджують у більш тонко градуйованих шкалах функціональних просторів, ніж класичні шкали Соболєва і Гельдера
(див., наприклад, \cite{DenkHieberPruess07,Lindemulder20,LosMikhailetsMurach21CPAA,LosMikhailetsMurach17CPAA,DongKim15,Hummel21,LeCronePrussWilke14}).
Використовують простори зі змішаною $L_p-L_q$ нормою, простори Лізоркіна-Трібеля, різні вагові простори, узагальнені простори Соболєва $H^{s,s/(2b);\varphi}$, де $2b$ -- параболічна вага. Для останніх теорія розв'язності скалярних параболічних задач побудована в основному в працях  \cite{LosMikhailetsMurach21CPAA,LosMikhailetsMurach17CPAA,LosMurach17OpenMath} і
викладена в монографії \cite{LosMikhailetsMurach21Monograph}. Випадок систем досліджено в \cite{DiachenkoLos22JEPE,Los17UMJ3}.

Основні результати цієї статті -- теореми про достатні умови приналежності узагальнених розв'язків досліджуваних у \cite{DiachenkoLos22JEPE} параболічних задач
до просторів $H^{s,s/2;\varphi}(\Omega)$ та умови неперервності розв'язків та вказаних їх частинних похідних, зокрема умови класичності узагальнених розв'язків.
Використання функціонального параметра $\varphi$ дозволяє отримати нові тонкі і точні умови гладкості розв'язків у порівнянні з класичними результатами
\cite{Solonnikov65,Ilin60, IlinKalashOleinik62}.
У скалярному випадку версії цих теорем встановлено в \cite{LosMikhailetsMurach21CPAA, Los16UMJ9, Los16UMJ11},
а для еліптичних крайових задач -- в \cite{MikhailetsMurach14,AnopDenkMurach21CPAA,MikhailetsMurach12BJMA2,AnopChepurukhinaMurac21Axioms}.
Окремий випадок параболічних крайових задач для систем з однорідними початковими даними Коші розглянуто в \cite{Los17UMJ3,Los20MFAT2}.

\textbf{2. Постановка задачі.}
Нехай довільно задані ціле
число $n\geq2$, дійсне число $\tau>0$
і обмежена область $G\subset\mathbb{R}^{n}$ з
нескінченно гладкою межею $\Gamma:=\partial G$.
Позначимо $\Omega:=G\times(0,\tau)$ --- відкритий циліндр в $\mathbb{R}^{n+1}$,
$S:=\Gamma\times(0,\tau)$~--- його бічна поверхня.
Тоді $\overline{\Omega}:=\overline{G}\times[0,\tau]$ і
$\overline{S}:=\Gamma\times[0,\tau]$ є замикання
$\Omega$ і $S$ відповідно. Будемо ототожнювати $G$ з нижньою основою $G\times\{0\}$ циліндра $\Omega$.

Розглянемо  таку параболічну початково-крайову задачу у циліндрі $\Omega$:
\begin{equation}\label{25f1}
\begin{split}
\partial_t u_j(x,t)+
\sum_{k=1}^{N}\sum_{|\alpha|\leq2}&a^{\alpha}_{j,k}(x,t)\,D^\alpha_x u_{k}(x,t)=f_{j}(x,t)\\
&\mbox{для всіх}\quad (x,t)\in\Omega\quad\mbox{та}\quad
j\in\{1,\dots,N\};
\end{split}
\end{equation}
\begin{equation}\label{25f2}
\begin{split}
\sum_{k=1}^{N}\sum_{|\alpha|\leq l_j}&b^{\alpha}_{j,k}(x,t)\,D^\alpha_x u_{k}(x,t)=g_{j}(x,t)\\
&\mbox{для всіх}\quad (x,t)\in S\quad\mbox{та}\quad
j\in\{1,\dots,N\};
\end{split}
\end{equation}
\begin{equation}\label{25f3}
u_{j}(x,t)\big|_{t=0}=h_{j}(x)\quad
\mbox{для всіх}\quad x\in G\quad\mbox{та}\quad j\in\{1,\ldots,N\}.
\end{equation}
Тут довільним чином вибрані натуральне число $N\geq2$ та числа $l_1,\ldots,l_N\in\{0,\,1\}$.
Всі коефіцієнти диференціальних виразів у формулах \eqref{25f1} та \eqref{25f2}
є нескінченно гладкими комплекснозначними функціями, заданими на $\overline{\Omega}$ і $\overline{S}$ відповідно; тобто усі
$a_{j,k}^{\alpha}\in C^{\infty}(\overline{\Omega})$ і
$b_{j,k}^{\alpha}\in C^{\infty}(\overline{S})$.
Використовуємо такі позначення
$D^\alpha_x:=D^{\alpha_1}_{1}\dots D^{\alpha_n}_{n}$, де $D_{k}:=i\,\partial/\partial{x_k}$ і $\partial_t:=\partial/\partial t$
для частинних похідних функцій, що залежать від $x=(x_1,\ldots,x_n)\in\mathbb{R}^{n}$ і $t\in\mathbb{R}$.
Тут $i$ -- уявна одиниця, $\alpha=(\alpha_1,...,\alpha_n)$ мультиіндекс, і $|\alpha|:=\alpha_1+\cdots+\alpha_n$.
У формулах \eqref{25f1} і \eqref{25f2} та їх аналогах підсумовування ведеться за цілими невід'ємними індексами
$\alpha_1,...,\alpha_n$, які задовольняють умову, вказану під знаком суми.

Припускаємо, що початково-крайова задача
\eqref{25f1}--\eqref{25f3} параболічна за Петровським
у циліндрі $\Omega$ (див. означення в \cite[розд.~1, \S~1]{Solonnikov65}).

Встановимо достатні умови глобальної та локальної регулярності узагальнених розв'язків задачі  \eqref{25f1}--\eqref{25f3} в термінах приналежності її правих частин узагальненим
анізотропним просторам Соболєва. Крім того отримаємо нові достатні умови класичності цих розв'язків.

Нагадаємо означення узагальненого анізотропного простору Соболєва на $\mathbb{R}^{k}$.
Через $\mathcal{M}$ позначимо клас усіх вимірних за Борелем функцій
$\varphi:[1,\infty)\rightarrow(0,\infty)$ таких, що:

(i) $\varphi$ і $1/\varphi$ обмежені на кожному відрізку $[1,c]$, де $1<c<\infty$;

(ii) $\varphi$ повільно змінна за Й.~Карамата на нескінченності, тобто
\begin{equation*}
\lim_{r\rightarrow\infty}\frac{\varphi(\lambda r)}{\varphi(r)}=1\quad\mbox{для кожного}\quad\lambda>0.
\end{equation*}

Нехай $s\in\mathbb{R}$ і $\varphi\in\mathcal{M}$.
За означенням, комплексний лінійний простір $H^{s,s/2;\varphi}(\mathbb{R}^{k})$, де $2\leq k\in\mathbb{Z}$, складається з усіх повільно зростаючих розподілів $w\in \mathcal{S}'(\mathbb{R}^{k})$ таких, що їх (повне) перетворення
Фур'є $\widetilde{w}$ є функцією, яка локально інтегровна на $\mathbb{R}^{k}$ за Лебегом і задовольняє умову
\begin{equation}\label{norm}
\begin{split}
&\|w\|_{H^{s,s/2;\varphi}(\mathbb{R}^{k})}:=\\
&\biggl(\;\int\limits_{\mathbb{R}^{k-1}}\int\limits_{\mathbb{R}}
\bigl(1+|\xi|^2+|\eta|\bigr)^{s}
\,\varphi^{2}\bigl((1+|\xi|^2+|\eta|)^{1/2}\bigr)\,
|\widetilde{w}(\xi,\eta)|^{2}\,d\xi\,d\eta\biggr)^{1/2}<\infty,
\end{split}
\end{equation}
де $\xi\in\mathbb{R}^{k-1}$ і $\eta\in\mathbb{R}$.
Цей простір гільбертовий і сепарабельний відносно норми \eqref{norm}.
Він є окремим випадком просторів $\mathcal{B}_{p,\mu}$, введених Л. Хермандером \cite[п.~2.2]{Hermander63}; а саме, $H^{s,s/2;\varphi}(\mathbb{R}^{k})=
\mathcal{B}_{p,\mu}$ за умови, що $p=2$ і
$$\mu(\xi,\eta)\equiv \bigl(1+|\xi|^2+|\eta|\bigr)^{s/2}\,\varphi\bigl((1+|\xi|^2+|\eta|)^{1/2}\bigr).
$$

Гільбертовий анізотропний простір $H^{s,s/2;\varphi}(\Omega)$ означається як простір звужень на $\Omega$ усіх розподілів з $H^{s,s/2;\varphi}(\mathbb{R}^{n+1})$, а
гільбертовий анізотропний простір $H^{s,s/2;\varphi}(S)$ на бічній поверхні циліндра означається за базовим простором $H^{s,s/2;\varphi}(\mathbb{R}^{n})$ за допомогою спеціальних
локальних карт на бічній поверхні циліндра (див. \cite[п.~1]{Los16JMathSci}).
Означення та основні властивості просторів $H^{s,s/2;\varphi}(W)$, де $W\in\{\Omega, S\}$, наведені, наприклад, в \cite[п.2]{DiachenkoLos22JEPE}.
Ізотропні простори $H^{s;\varphi}(G)$ і $H^{s;\varphi}(\Gamma)$, задані на основі $G$ циліндра та лінії $\Gamma$ з'єднання основи і бічної поверхні відповідно, означено в
\cite[п.2.1.1, 3.2.1]{MikhailetsMurach14}, \cite{MikhailetsMurach12BJMA2}.
Якщо $\varphi(\cdot)=1$, то ці простори є соболєвськими.
В цьому випадку прибираємо індекс $\varphi$ у їх позначеннях.

Результати цієї роботи спираються на теорему про ізоморфізми \cite[теорема 4.1]{DiachenkoLos22JEPE} для задачі \eqref{25f1}--\eqref{25f3}.
Сформулюємо для зручності цей результат.

Покладемо
\begin{equation}\label{25f4}
A_{j,k}(x,t,D_x,\partial_t):=\delta_{j,k}\partial_t+
\sum_{|\alpha|\leq 2}a^{\alpha}_{j,k}(x,t)\,D^\alpha_x
\end{equation}
та
\begin{equation}\label{25f5}
B_{j,k}(x,t,D_x):=\sum_{|\alpha|\leq l_{j}}b^{\alpha}_{j,k}(x,t)\,D^\alpha_x
\end{equation}
для всіх $j,k\in\{1,\ldots,N\}$. Тут $\delta_{j,k}$ -- символ Кронекера.
Скориставшись позначеннями \eqref{25f4} та \eqref{25f5} запишемо
всі рівності в \eqref{25f1} та \eqref{25f2} у такій еквівалентній формі:
\begin{equation*}
\sum_{k=1}^{N}A_{j,k}(x,t,D_x,\partial_t)u_{k}(x,t)=f_{j}(x,t)
\end{equation*}
та
\begin{equation*}
\sum_{k=1}^{N}B_{j,k}(x,t,D_x)u_{k}(x,t)=
g_{j}(x,t).
\end{equation*}
Покладемо $u:=(u_1,\ldots,u_N)$, $f:=(f_1,\ldots,f_N)$, $g:=(g_1,\ldots,g_N)$ та $h:=(h_1,\ldots,h_N)$.
Запишемо систему \eqref{25f1} та граничні умови \eqref{25f2} в матричній формі
$Au=f$ та $Bu=g$; тут
$$
A:=(A_{j,k}(x,t,D_x,\partial_t))_{j,k=1}^N\quad\mbox{та}\quad
B:=\bigl(B_{j,k}(x,t,D_x)\bigr)_{j,k=1}^N
$$
-- матричні диференціальні оператори. Пов'яжемо із задачею \eqref{25f1}--\eqref{25f3} лінійне відображення
\begin{equation}\label{25f4.1}
u\mapsto\Lambda
u:=\bigl(Au,Bu,u|_{t=0}\bigr),\quad\mbox{де}\;\,u\in
\bigl(C^{\infty}(\overline{\Omega})\bigr)^{N}.
\end{equation}

\textbf{Твердження 1} \cite[теорема 4.1]{DiachenkoLos22JEPE}.
\it Нехай довільно задані параметри: числовий
$s\geq2$ і функціональний $\varphi\in\nobreak\mathcal{M}$.
У випадку $s=2$ додатково припустимо, що
$\varphi$ є зростаючою (в нестрогому сенсі) функцією.
Тоді відображення \eqref{25f4.1} продовжується єдиним чином (за неперервністю) до ізоморфізму
\begin{equation}\label{25f11}
\Lambda :\,\bigl(H^{s,s/2;\varphi}(\Omega)\bigr)^N\leftrightarrow
\mathcal{Q}^{s-2,s/2-1;\varphi}.
\end{equation}
\rm
Тут $\mathcal{Q}^{s-2,s/2-1;\varphi}$ гільбертовий простір правих частин задачі (див. \cite[п.4]{DiachenkoLos22JEPE}).
А саме, для $s\notin E$, де
\begin{equation}\label{setE}
E:=\{2l+l_j+3/2:j,l\in\mathbb{Z},\;1\leq j\leq N,\;l\geq0\}
\cap(2,\infty),
\end{equation}
він утворений такими векторами $(f,g,h)$ з простору
\begin{align*}
\mathcal{H}^{s-2,s/2-1;\varphi}:=
&\bigl(H^{s-2,s/2-1;\varphi}(\Omega)\bigr)^N\\
&\oplus\bigoplus_{j=1}^{N}H^{s-l_j-1/2,(s-l_j-1/2)/2;\varphi}(S)
\oplus\bigl(H^{s-1;\varphi}(G)\bigr)^N,
\end{align*}
що задовольняють природні умови узгодження правих частин параболічної задачі \eqref{25f1}--\eqref{25f3}.
У випадку $s\in E$ гільбертовий простір $\mathcal{Q}^{s-2,s/2-1;\varphi}$ означається за допомогою квадратичної інтерполяції з числовим параметром $1/2$:
\begin{equation}\label{16f10}
\begin{split}
&\mathcal{Q}^{s-2,s/2-1;\varphi}:= \\
&:=
\bigl[\mathcal{Q}^{s-\varepsilon-2,(s-\varepsilon)/2-1;\varphi},
\mathcal{Q}^{s+\varepsilon-2,(s+\varepsilon)/2-1;\varphi}\bigr]_{1/2}.
\end{split}
\end{equation}
Тут $\varepsilon\in(0,1/2)$ -- довільне число. Означений у такий спосіб простір не залежить з точністю до еквівалентності норм від вибору числа $\varepsilon$.

\textbf{3. Основні результати.}
Це -- достатні умови глобальної та локальної регулярності узагальненого розв'язку параболічної задачі \eqref{25f1}--\eqref{25f3} в узагальнених просторах
Соболєва, а також -- достатня умова класичності цього розв'язку.
Дамо його означення.

Нехай усі компоненти правих частин $f$, $g$ та $h$ задачі є довільними розподілами на $\Omega$, $S$ і $G$ відповідно.
Вектор-функцію
$u\in \bigl(H^{2,1}(\Omega)\bigr)^N$ називаємо узагальненим розв'язком задачі \eqref{25f1}--\eqref{25f3}, якщо
\begin{equation}\label{Lambdau=f-Murach}
\Lambda u=(f,g,h),
\end{equation}
де $\Lambda$~--- обмежений оператор \eqref{25f11} для $s=2$ і $\varphi=1$. З рівності \eqref{Lambdau=f-Murach} випливає, що
\begin{equation}\label{16f12}
(f,g,h)\in\mathcal{Q}^{0,0}.
\end{equation}
З ізоморфізму \eqref{25f11} випливає (див. також \cite[теорема~5.7]{ZhitarashuEidelman98}), що параболічна задача \eqref{25f1}--\eqref{25f3} має єдиний узагальнений розв'язок  $u\in \bigl(H^{2,1}(\Omega)\bigr)^N$ для кожного вектора \eqref{16f12}.

\textbf{Теорема 1.} \it
Припустимо, що вектор-функція $u\in \bigl(H^{2,1}(\Omega)\bigr)^N$
є узагальненим розв'язком параболічної задачі
\eqref{25f1}--\eqref{25f3}, праві частини якої задовольняють умову
$$
(f,g,h)\in
\mathcal{Q}^{s-2,s/2-1;\varphi}
$$
для деяких $s\geq2$ і $\varphi\in\mathcal{M}$ (у випадку $s=2$ додатково припускаємо, що
$\varphi$ є зростаючою (в нестрогому сенсі) функцією). Тоді $u\in \bigl(H^{s,s/2;\varphi}(\Omega)\bigr)^N$.
\rm

Цей результат про достатню умову глобальної (тобто в усьому циліндрі $\Omega$ аж до його межі) регулярності розв'язку є прямим наслідком твердження~1.

Тепер сформулюємо локальну версію цієї теореми.
Нехай $U$~--- відкрита множина в $\mathbb{R}^{n+1}$ така, що $\nobreak{\Omega_0:=U\cap\Omega\neq\varnothing}$ і $U\cap\Gamma=\varnothing$.
Покладемо $\Omega':=U\cap\partial\overline{\Omega}$, $S_0:=U\cap S$, $S':=U\cap \{(x,\tau):x\in\Gamma\}$ і $G_0:=U\cap G$.
Позначимо через $H^{s,s/2;\varphi}_{\mathrm{loc}}(\Omega_0,\Omega')$ лінійний простір усіх розподілів $u$ на $\Omega$ таких, що $\chi u\in H^{s,s/2;\varphi}(\Omega)$ для кожної функції $\chi\in C^\infty (\overline\Omega)$, яка задовольняє умову $\mathrm{supp}\,\chi\subset\Omega_0\cup\Omega'$.
Аналогічно, позначимо через $H^{s,s/2;\varphi}_{\mathrm{loc}}(S_0,S')$ лінійний простір усіх розподілів $v$ на $S$ таких, що $\chi v\in H^{s,s/2;\varphi}(S)$ для будь-якої функції $\chi\in C^\infty (\overline S)$, яка задовольняє умову $\mathrm{supp}\,\chi\subset S_0\cup S'$. Нарешті, $H^{s;\varphi}_{\mathrm{loc}}( G_0)$ позначає лінійний простір усіх розподілів $w$ на $G$ таких, що $\chi w\in H^{s;\varphi}(G)$ для кожної функції $\chi\in C^\infty (\overline G)$, яка задовольняє умову $\mathrm{supp}\,\chi\subset G_0$
(див., наприклад, \cite[пункт 4]{LosMikhailetsMurach21CPAA}).

\textbf{Теорема 2.} \it
Нехай $s\geq2$ і $\varphi\in\mathcal{M}$ (у випадку $s=2$ додатково припускаємо, що
$\varphi$ є зростаючою (в нестрогому сенсі) функцією). Припустимо, що вектор-функція $u\in \bigl(H^{2,1}(\Omega)\bigr)^N$
є узагальненим розв'язком параболічної задачі \eqref{25f1}--\eqref{25f3}, праві частини якої задовольняють такі умови:
\begin{equation}\label{16f13}
f\in\bigl(H^{s-2,s/2-1;\varphi}_{\mathrm{loc}}(\Omega_0,\Omega')\bigr)^N,
\end{equation}
\begin{equation}\label{16f14}
g\in\bigoplus_{j=1}^{N}H^{s-l_j-1/2,(s-l_j-1/2)/2;\varphi}_{\mathrm{loc}}(S_0,S'),
\end{equation}
\begin{equation}\label{16f15}
h\in\bigl(H^{s-1;\varphi}_{\mathrm{loc}}(G_0)\bigr)^N.
\end{equation}
Тоді $u\in\bigl(H^{s,s/2;\varphi}_{\mathrm{loc}}(\Omega_0,\Omega')\bigr)^N$.
\rm

Якщо $\Omega'=\varnothing$, то теорема~2
стверджує, що регулярність узагальненого розв'язку підвищується в околах внутрішніх точок замкненого циліндра $\overline{\Omega}$.
Якщо $G_0=\varnothing$, то ця теорема стверджує, що регулярність розв'язку $u(x,t)$ підвищується при $t>0$.
В цих випадках вона є наслідком \cite[теорема~2]{Los17UMJ3}.

Зазначимо, що припущення $U\cap\Gamma=\varnothing$ істотне, оскільки без нього висновок теореми~2 є взагалі хибним.
Для його істинності у цьому випадку, треба накласти на праві частини задачі \eqref{25f1}--\eqref{25f3} на множині $U\cap\Gamma$ деякі додаткові умови узгодження.

Узагальнені простори Соболєва дозволяють отримати більш тонкі, ніж у випадку класичних просторів Соболєва, достатні умови неперервності
узагальненого розв'язку $u$ та його узагальнених похідних заданого порядку на множині $\Omega_0\cup\Omega'$.

Подібно до \cite[с.3617]{LosMikhailetsMurach21CPAA} (див. також \cite[зауваження~2.1]{LosMikhailetsMurach21Monograph})
узагальнену функцію $v\in\mathcal{D}'(\Omega)$ називаємо неперервною на множині $\Omega_0\cup\Omega'$, якщо існує неперервна функція $v_{0}$ на $\Omega_0\cup\Omega'$ така, що
\begin{equation}\label{rem-to-th4.4}
v(\omega)=\int\limits_{\Omega_0}v_{0}(x,t)\,\omega(x,t)\,dxdt
\end{equation}
для довільної функції $\omega\in C^{\infty}(\Omega)$, носій якої задовольняє умову $\mathrm{supp}\,\omega\subset\Omega_0$.
Тут $v(\omega)$~--- значення функціонала $v$ на функції~$\omega$.

\textbf{Теорема 3.} \it
Нехай задане довільне ціле число $p\geq0$.
Припустимо, що вектор-функція $u\in \bigl(H^{2,1}(\Omega)\bigr)^N$ є узагальненим розв'язком параболічної задачі \eqref{25f1}--\eqref{25f3}, праві частини якої задовольняють умови \eqref{16f13}--\eqref{16f15} для $s:=p+1+n/2$ і деякого функціонального параметра $\varphi\in\mathcal{M}$,
що задовольняє умову
\begin{equation}\label{9f4.7}
\int\limits_{1}^{\,\infty}\;\frac{dr}{r\,\varphi^2(r)}<\infty,
\end{equation}
причому у випадку $s=2$ додатково припускаємо, що функція $\varphi$ зростає (в нестрогому сенсі).
Тоді розв'язок $u(x,t)=(u_1(x,t),\ldots,u_N(x,t))$ і кожна його узагальнена частинна похідна
$D_{x}^{\alpha}\partial_{t}^{\beta}u(x,t)=(D_{x}^{\alpha}\partial_{t}^{\beta}u_1(x,t),\ldots,D_{x}^{\alpha}\partial_{t}^{\beta}u_N(x,t))$,
де $|\alpha|+2\beta\leq p$,  неперервні на множині $\Omega_0\cup\Omega'$.
\rm

Для цієї теореми правильні версії зауважень~1 і~2 з \cite{Los17UMJ3}.
А саме, умова \eqref{9f4.7} в теоремі є точною.
Крім того, якщо переформулювати теорему~3 на випадок класичних просторів Соболєва ($\varphi=1$),
то \eqref{9f4.7} не виконується і  доведеться  замінити умови \eqref{16f13}--\eqref{16f15} для $s:=p+1+n/2$ на більш сильні.
Потрібно стверджувати, що ці умови виконуються для деякого $s>p+1+n/2$.

Теорема~3 дозволяє отримати нові й тонкі достатні умови класичності узагальненого розв'язку задачі \eqref{25f1}--\eqref{25f3}.
Сформулюємо означення класичного розв'язку цієї задачі.

Нехай $l_0:=\max\{l_1,\dots,l_N\}$. Покладемо
\begin{equation*}
S_{\varepsilon}:=\{x\in\Omega:\mbox{dist}(x,S)<\varepsilon\},\quad
G_{\varepsilon}:=\{x\in\Omega:\mbox{dist}(x,G)<\varepsilon\},
\end{equation*}
де число $\varepsilon>0$.

Узагальнений розв'язок $u\in \bigl(H^{2,1}(\Omega)\bigr)^N$ задачі \eqref{25f1}--\eqref{25f3} називаємо класичним, якщо узагальнені частинні похідні вектор-функції $u=u(x,t)$ задовольняють такі три умови:
\begin{itemize}
\item [(a)] $D^\alpha_x\partial^\beta_t u$ неперервна на $\Omega$, якщо $0\leq|\alpha|+2\beta\leq 2$;
\item [(b)] $D^\alpha_x u$ неперервна на $S_{\varepsilon}\cup S$ для деякого числа $\varepsilon>0$, якщо $0\leq|\alpha|\leq l_0$;
\item [(c)] $u$ неперервна на $G_{\varepsilon}\cup G$ для деякого числа $\varepsilon>0$.
\end{itemize}

Якщо розв'язок $u=u(x,t)$ задачі \eqref{25f1}--\eqref{25f3} класичний, то її ліві частини є неперервними функціями на відповідних множинах.

\textbf{Теорема 4.} \it
Припустимо, що вектор-функція $u\in \bigl(H^{2,1}(\Omega)\bigr)^N$ є узагальненим розв'язком параболічної задачі \eqref{25f1}--\eqref{25f3}, праві частини якої задовольняють такі умови:
\begin{equation}\label{24f12}
\begin{aligned}
f\in\, &\bigl(H_{\mathrm{loc}}^{1+n/2,\,1/2+n/4;\varphi}
(\Omega,\varnothing)\bigr)^N\cap\\
&\cap
\bigl(H_{\mathrm{loc}}^{l_0-1+n/2,\,l_0/2-1/2+n/4;\varphi}
(S_{\varepsilon},S)\bigr)^N\cap\\
&\cap
\bigl(H_{\mathrm{loc}}^{-1+n/2,\,-1/2+n/4;\varphi}
(G_{\varepsilon},G)\bigr)^N,
\end{aligned}
\end{equation}
\begin{equation}\label{24f12bis}
g\in\bigoplus_{j=1}^{N}H_{\mathrm{loc}}^{l_0+n/2-l_j+1/2,\,l_0/2+n/4-l_j/2+1/4;\varphi}
(S,\varnothing),
\end{equation}
\begin{equation}\label{18f9}
h\in\bigl(H_{\mathrm{loc}}^{n/2;\varphi}(G)\bigr)^N,
\end{equation}
з деяким функціональним параметром $\varphi\in\mathcal{M}$, який задовольняє умови теореми~3. Тоді розв'язок $u$ класичний.
\rm

\textbf{6. Доведення.}

У скалярному випадку $N=1$ теореми~2, 3 і~4 встановлено в \cite[розд.~6]{LosMikhailetsMurach21CPAA} (див. також монографію \cite[розд.~3.4, 3.5]{LosMikhailetsMurach21Monograph}).
Ми узагальнюємо розроблені там методи доведення на випадок систем.

\textbf{Доведення теореми 2.}
Першим кроком доведемо, що з умов \eqref{16f13}--\eqref{16f15}  випливає істинність імплікації
\begin{equation}\label{16f54}
\begin{aligned}
&u\in \bigl(H^{s-\lambda,(s-\lambda)/2;\varphi}_{\mathrm{loc}}
(\Omega_0,\Omega')\bigr)^N\Longrightarrow\\
&\Longrightarrow u\in\bigl(H^{s-\lambda+1,(s-\lambda+1)/2;\varphi}_{\mathrm{loc}}
(\Omega_0,\Omega')\bigr)^N
\end{aligned}
\end{equation}
для всіх цілих $\lambda\geq1$, що задовольняють умову $s-\lambda+1>2$.

Нехай  $\chi\in C^\infty(\overline\Omega)$ є довільною функцією такою, що  $\mbox{supp}\,\chi\subset\Omega_0\cup\Omega'$. Для $\chi$ існує
функція $\eta\in C^\infty(\overline\Omega)$ така, що $\mbox{supp}\,\eta\subset\Omega_0\cup\Omega'$ і
$\eta=1$ в околі $\mbox{supp}\,\chi$.
Переставивши диференціальні  оператори $A$ і $B$  з оператором множення на функцію $\chi$, отримаємо
\begin{equation}\label{16f55}
\begin{aligned}
\Lambda(\chi u)&=\Lambda(\chi\eta u)=
\chi\Lambda(\eta u)+ \Lambda'(\eta u)=\\
&=\chi\Lambda u+\Lambda'(\eta u)=\chi(f,g,h)+
\Lambda'(\eta u).
\end{aligned}
\end{equation}
Тут $\Lambda':=(A',B',0)$ -- оператор, елементи $A'$ і $B'$ якого мають вигляд
$$
A':=\biggl(\sum_{|\alpha|\leq 1}
a^{\alpha}_{j,k,1}(x,t)\,D^\alpha_x\biggr)_{j,k=1}^N\quad\mbox{і}\quad
B':=\bigl(b_{j,k,1}(x,t)\bigr)_{j,k=1}^N,
$$
де $a^{\alpha}_{j,k,1}\in C^{\infty}(\overline{\Omega})$ і $b_{j,k,1}\in C^{\infty}(\overline{S})$.
Іншими словами, $A'$ і $B'$~--- деякі матричні диференціальні оператори, порядки кожної компоненти яких принаймні на одиницю менші, ніж порядки відповідних
компонент операторів $A$ і $B$.

Для кожного  $\sigma\geq1$ оператор $\Lambda'$ неперервно діє на парі просторів
\begin{equation}\label{16f59}
\Lambda':\bigl(H^{\sigma,\sigma/2;\varphi}(\Omega)\bigr)^N\rightarrow
\mathcal{H}^{\sigma-1,\,\sigma/2-1/2;\,\varphi}.
\end{equation}
У випадку $\varphi(\cdot)\equiv1$ цей факт відомий. Він випливає з властивостей операторів диференціювання та операторів сліду у просторах Соболєва.
Звідси неперервність оператора \eqref{16f59} у загальному випадку отримується методом квадратичної інтерполяції.
(Її викладено, наприклад, у монографії \cite[п.1.1]{MikhailetsMurach14}).
Для $\sigma>1$ і довільного $\varphi\in\mathcal{M}$ це випливає з інтерполяційних формул
\cite[Теорема~2]{Los16JMathSci}, \cite[Лема~2]{Los16JMathSci} (з $\Omega$ замість $\Pi$) та \cite[Теореми~1.5, 1.14(i), 3.2]{MikhailetsMurach14}.
Для $\sigma=1$ і зростаючої функції $\varphi\in\mathcal{M}$ --- з \cite[Леми~2 і 3]{Los16UMJ6}, \cite[Теорема~4.1]{MikhailetsMurach15ResMath1}
та \cite[Теорема~1.5]{MikhailetsMurach14}.

З умов \eqref{16f13}--\eqref{16f15} теореми випливає, що
\begin{equation}\label{16f59a}
\chi\,(f,g,h)\in\mathcal{H}^{s-2,s/2-1;\varphi}.
\end{equation}
Врахувавши неперервність оператора $\Lambda'$ на парі просторів \eqref{16f59}, де $\sigma:=s-\lambda$, маємо
\begin{align*}
&u\in \bigl(H^{s-\lambda,(s-\lambda)/2;\varphi}_{\mathrm{loc}}
(\Omega_0,\Omega')\bigr)^N \Longrightarrow\\
&\Longrightarrow\Lambda'(\eta u)\in
\mathcal{H}^{s-\lambda-1,(s-\lambda-1)/2;\varphi}.
\end{align*}
З цієї імплікації, включення \eqref{16f59a} та формули \eqref{16f55} випливає, що
\begin{equation}\label{16f60}
\begin{aligned}
&u\in \bigl(H^{s-\lambda,(s-\lambda)/2;\varphi}_{\mathrm{loc}}
(\Omega_0,\Omega')\bigr)^N\Longrightarrow\\
&\Longrightarrow\Lambda(\chi u)\in
\mathcal{H}^{s-\lambda-1,(s-\lambda-1)/2;\varphi}.
\end{aligned}
\end{equation}

Далі покажемо, що з умов $U\cap\Gamma=\varnothing$ і $\Lambda(\chi u)\in
\mathcal{H}^{\sigma-2,\sigma/2-1;\varphi}$ випливає включення
\begin{equation}\label{16f61}
\Lambda(\chi u)\in\mathcal{Q}^{\sigma-2,\sigma/2-1;\varphi}
\end{equation}
для будь-якого  $\sigma\geq2$. Цей факт дасть можливість довести \eqref{16f54}, скориставшись теоремою~1.

Оскільки $U\cap\Gamma=\varnothing$, то $\mathrm{dist}(\mathrm{supp}\,\chi,\Gamma)>0$. Це означає, що $\Lambda(\chi u)=0$ в деякому околі множини $\Gamma$.
Тому вектор-функція $\Lambda(\chi u)$ задовольняє умови узгодження правих частин параболічної задачі \eqref{25f1}--\eqref{25f3} (див., наприклад, \cite[п.4]{DiachenkoLos22JEPE}).
Отже, згідно з означенням простору $\mathcal{Q}^{\sigma-2,\sigma/2-1;\varphi}$ у випадку $\sigma\notin E$ виконується включення
$\Lambda(\chi u)\in\mathcal{Q}^{\sigma-2,\sigma/2-1;\varphi}$.

Розглянемо випадок $\sigma\in E$.
Нехай $\chi_1\in C^{\infty}(\overline\Omega)$ така, що $\chi_1=0$ в околі множини $\Gamma$ і  $\chi_1=1$ в околі $\mathrm{supp}\,\chi$.
З попередніх міркувань випливає, що відображення $M_{\chi_1}:F\mapsto\chi_1 F$ є обмеженим оператором на просторах
\begin{equation}\label{16f63}
\begin{aligned}
M_{\chi_1}&:
\mathcal{H}^{\sigma\pm\varepsilon-2,(\sigma\pm\varepsilon)/2-1;\varphi}
\to\\
&\to\mathcal{Q}^
{\sigma\pm\varepsilon-2,(\sigma\pm\varepsilon)/2-1;\varphi},
\end{aligned}
\end{equation}
якщо $0<\varepsilon<1/2$, оскільки $\sigma\pm\varepsilon\notin E$.
Інтерполюючи з числовим параметром $1/2$  оператори \eqref{16f63} маємо обмежений оператор
\begin{equation*}
\begin{aligned}
M_{\chi_1}&:
\bigl[
\mathcal{H}^{\sigma-\varepsilon-2,(\sigma-\varepsilon)/2-1;\varphi},
\mathcal{H}^{\sigma+\varepsilon-2,(\sigma+\varepsilon)/2-1;\varphi}
\bigr]_{1/2}\to
\\
&\to\bigl[
\mathcal{Q}^{\sigma-\varepsilon-2,(\sigma-\varepsilon)/2-1;\varphi},
\mathcal{Q}^{\sigma+\varepsilon-2m,(\sigma+\varepsilon)/2-1;\varphi}
\bigr]_{1/2}.
\end{aligned}
\end{equation*}
Згідно з інтерполяційними формулами для просторів $\mathcal{H}^{\sigma\pm\varepsilon-2,(\sigma\pm\varepsilon)/2-1;\varphi}$ (див. \cite[Лема 6.4]{LosMikhailetsMurach21CPAA})
та \eqref{16f10} цей оператор діє на парі просторів
\begin{equation}\label{16f65}
M_{\chi_1}:\mathcal{H}^{\sigma-2,\sigma/2-1;\varphi}\to
\mathcal{Q}^{\sigma-2,\sigma/2-1;\varphi}.
\end{equation}
Оскільки $\chi_1$ така, що $\chi_1=1$ в околі $\mathrm{supp}\,\chi$, то $\chi_1\Lambda(\chi u)=\Lambda(\chi u)$.
Нагадаємо, що $\Lambda(\chi u)\in\mathcal{H}^{\sigma-2,\sigma/2-1;\varphi}$.
Тоді на підставі \eqref{16f65} виконується включення
$$
\Lambda(\chi u)=\chi_1\Lambda(\chi u)\in \mathcal{Q}^{\sigma-2,\sigma/2-1;\varphi}.
$$
Випадок $\sigma\in E$ розглянуто.

З формул \eqref{16f60} і \eqref{16f61}, де $\sigma=s-\lambda+1$, випливає за теоремою~1, що
\begin{align*}
&u\in\bigl(H^{s-\lambda,(s-\lambda)/2;\varphi}_{\mathrm{loc}}
(\Omega_0,\Omega')\bigr)^N\Longrightarrow\\
&\Longrightarrow\Lambda(\chi u)\in
\mathcal{H}^{s-\lambda-1,(s-\lambda-1)/2;\varphi}\Longrightarrow\\
&\Longrightarrow \Lambda(\chi u)\in
\mathcal{Q}^{s-\lambda-1,(s-\lambda-1)/2;\varphi}\Longrightarrow\\
&\Longrightarrow
\chi u\in\bigl(H^{s-\lambda+1,(s-\lambda+1)/2;\varphi}(\Omega)\bigr)^N.
\end{align*}
Оскільки $\chi\in C^\infty(\overline\Omega)$ є довільною функцією, підпорядкованою умові $\mbox{supp}\,\chi\subset\Omega_0\cup\Omega'$, то
останнє включення означає, що $u\in\bigl(H^{s-\lambda+1,(s-\lambda+1)/2;\varphi}_{\mathrm{loc}}
(\Omega_0,\Omega')\bigr)^N$.
Оскільки за умовою теореми $\chi u\in\bigl(H^{2,1}(\Omega)\bigr)^N$ й, крім того, $s-\lambda+1>2$, то можна застосувати теорему~1.
Імплікацію \eqref{16f54} доведено.

За її допомогою доведемо потрібне включення $u\in\nobreak\bigl(H^{s,s/2;\varphi}_{\mathrm{loc}}(\Omega_0,\Omega')\bigr)^N$.
Нагадаємо, що $s\geq2$. Нехай спочатку $s\notin\mathbb{Z}$. Тоді існує ціле число $\lambda_{0}$ таке, що
\begin{equation*}\label{16f66}
s-\lambda_{0}<2<s-\lambda_{0}+1.
\end{equation*}
За умовою теореми $u\in\bigl(H^{2,1}(\Omega)\bigr)^N$. Використавши \eqref{16f54} послідовно для значень
$\lambda:=\nobreak\lambda_{0}$, $\lambda:=\lambda_{0}-1$,..., $\lambda:=1$, робимо висновок, що
\begin{align*}
&u\in\bigl(H^{2,1}(\Omega)\bigr)^N\subset
\bigl(H^{s-\lambda_{0},(s-\lambda_{0})/2;\varphi}_
{\mathrm{loc}}(\Omega_0,\Omega')\bigr)^N\Longrightarrow\\
&\Longrightarrow u\in
\bigl(H^{s-\lambda_{0}+1,(s-\lambda_{0}+1)/2;\varphi}_
{\mathrm{loc}}(\Omega_0,\Omega')\bigr)^N
\Longrightarrow\ldots\\
&\ldots\Longrightarrow
u\in\bigl(H^{s,s/2;\varphi}_{\mathrm{loc}}(\Omega_0,\Omega')\bigr)^N.
\end{align*}
Отже, потрібне включення доведено у розглянутому випадку.

Перейдемо до випадку коли $s\in\mathbb{Z}$ і $s>2$.
Використаємо щойно отриманий результат. Нехай $0<\varepsilon\ll1$. Тоді
$\nobreak{s-\varepsilon\notin\mathbb{Z}}$ і $s-\varepsilon>2$. Згідно з цим результатом
$$
u\in\bigl(H^{s-\varepsilon,(s-\varepsilon)/2;\varphi}_
{\mathrm{loc}}(\Omega_0,\Omega')\bigr)^N.
$$
Застосувавши імплікацію \eqref{16f54}, де $\lambda=1$, робимо висновок, що
\begin{align*}
&u\in\bigl(H^{s-\varepsilon,(s-\varepsilon)/2;\varphi}_
{\mathrm{loc}}(\Omega_0,\Omega')\bigr)^N\subset
\bigl(H^{s-1,(s-1)/2;\varphi}_{\mathrm{loc}}(\Omega_0,\Omega')\bigr)^N
\Longrightarrow\\
&\Longrightarrow u\in\bigl(H^{s,s/2;\varphi}_{\mathrm{loc}}(\Omega_0,\Omega')\bigr)^N.
\end{align*}

Випадок $s=2$ розглянемо окремо, оскільки в ньому є таке додаткове припущення: функція $\varphi\in\mathcal{M}$  зростає.
З умов \eqref{16f13}--\eqref{16f15} теореми випливає включення
\begin{equation}\label{16f59a1}
\chi\,(f,g,h)\in\mathcal{H}^{0,0;\varphi}.
\end{equation}
для довільної функції $\chi\in C^\infty(\overline\Omega)$, яка задовольняє умову $\mbox{supp}\,\chi\subset\Omega_0\cup\Omega'$.
З неперервності оператора \eqref{16f59}, де $\sigma=1$, випливає імплікація
\begin{equation*}
u\in \bigl(H^{1,1/2;\varphi}_{\mathrm{loc}}
(\Omega_0,\Omega')\bigr)^N \Longrightarrow
\Lambda'(\eta u)\in\mathcal{H}^{0,0;\varphi},
\end{equation*}
де функція $\eta$ така, як і раніше у доведенні.
На підставі умови $u\in\bigl(H^{2,1}(\Omega)\bigr)^N$, формули \eqref{16f55}, включення \eqref{16f59a1} та цієї імплікації робимо висновок, що
\begin{equation}\label{16f59a2}
\begin{aligned}
&u\in\bigl(H^{2,1}(\Omega)\bigr)^N\subset\bigl(H^{1,1/2;\varphi}_
{\mathrm{loc}}(\Omega_0,\Omega')\bigr)^N\Longrightarrow\\
&\Longrightarrow\Lambda(\chi u)\in
\mathcal{H}^{0,0;\varphi}.
\end{aligned}
\end{equation}
З формул \eqref{16f59a2} і \eqref{16f61}, де $\sigma=2$, та теореми~1 випливає, що
\begin{align*}
&u\in\bigl(H^{2,1}(\Omega)\bigr)^N\Longrightarrow\Lambda(\chi u)\in
\mathcal{H}^{0,0;\varphi}\Longrightarrow\\
&\Longrightarrow \Lambda(\chi u)\in
\mathcal{Q}^{0,0;\varphi}\Longrightarrow
\chi u\in\bigl(H^{2,1;\varphi}(\Omega)\bigr)^N.
\end{align*}

Теорему 2 доведено.

\textbf{Доведення теореми 3.}
Нагадаємо, що $n\geq2$, тому $p+1+n/2\geq2$.
Згідно з теоремою~2 виконується включення $u\in\bigl(H^{s,s/2;\varphi}_{\mathrm{loc}}(\Omega_0,\Omega')\bigr)^N$, де  $s=p+1+n/2$ і $\varphi$ задовольняє \eqref{9f4.7}.
У праці \cite[доведення теореми~4.3]{LosMikhailetsMurach21CPAA} (див. також \cite[доведення теореми~2.6]{LosMikhailetsMurach21Monograph})
встановлено такий результат: якщо узагальнена функція $v$ належить простору $H^{s,s/2;\varphi}_{\mathrm{loc}}(\Omega_0,\Omega')$, де $s=p+1+n/2$, а $\varphi$ задовольняє \eqref{9f4.7}, то вона разом з усіма її узагальненими частинними похідними $D_{x}^{\alpha}\partial_{t}^{\beta}v(x,t)$, де $|\alpha|+2\beta\leq p$,\, є неперервною на множині $\Omega_{0}\cup\Omega'_{0}$. Звідси і випливає висновок теореми~3.

\textbf{Доведення теореми 4.} Треба показати, що $u$ задовольняє умови (a)--(c) означення класичного розв'язку.
З умови~\eqref{24f12}, а саме, із включення
$$
f\in\,\bigl(H_{\mathrm{loc}}^{1+n/2,\,1/2+n/4;\varphi}
(\Omega,\varnothing)\bigr)^N
$$
випливає на підставі теореми~3 у випадку, коли $\Omega_0=\Omega$, $\Omega'=S_0=S'=G_0=\varnothing$ і $p=2$, що $u$ задовольняє умову (a).

Включення
$$
f\in \bigl(H_{\mathrm{loc}}^{l_0-1+n/2,\,l_0/2-1/2+n/4;\varphi}
(S_{\varepsilon},S)\bigr)^N
$$
і умова \eqref{24f12bis} тягнуть за собою виконання умови (b) для $u$ з огляду на теорему~3 у випадку,
коли $\Omega_0=S_\varepsilon$, $\Omega'=S_0=S$, $S'=G_0=\varnothing$ і $p=l_0$.

Нарешті, $u$ задовольняє умову (c) на підставі включення
$$
f\in\,\bigl(H_{\mathrm{loc}}^{-1+n/2,\,-1/2+n/4;\varphi}
(G_{\varepsilon},G)\bigr)^N
$$
і умови \eqref{18f9} згідно з теоремою~3 у випадку, коли $\Omega_0=G_{\varepsilon}$, $\Omega'=G_0=G$, $S_0=S'=\varnothing$ і $p=0$.

Теорему 4 доведено.

\end{document}